\documentclass[12pt,leqno]{article}
\usepackage{amsfonts}
\usepackage{amsmath, amsthm, amsfonts, amssymb, color}
\usepackage{mathrsfs}
\usepackage{color}
\theoremstyle{definition}
\newtheorem{defn}{Definition}[section]
\newcommand{\scr}[1]{\mathscr #1}
\definecolor{wco}{rgb}{0.5,0.2,0.3}

\numberwithin{equation}{section} \theoremstyle{remark}

\newcommand{\ua}{\uparrow}

\title{{\bf Bismut   Formula for  Intrinsic Derivative  of DDSDEs with Singular Interactions}\footnote{The author is supported by  NNSFC(12301180),
RGC (21301925), NSFC/RGC JRS N-CityU165/25  and Research Centre for Nonlinear Analysis at Hong Kong PolyU.} }
\author{
{\bf Panpan Ren}\\
Mathematics Department,  City University of Hong Kong,  Hong Kong,  China\\
panparen@cityu.edu.hk
 }
\begin{document}
\allowdisplaybreaks
\def\R{\mathbb R}  \def\ff{\frac} \def\ss{\sqrt} \def\B{\mathbf
B}
\def\N{\mathbb N} \def\kk{\kappa} \def\m{{\bf m}}
\def\ee{\varepsilon}\def\ddd{D^*}
\def\dd{\delta} \def\DD{\Delta} \def\vv{\varepsilon} \def\rr{\rho}
\def\<{\langle} \def\>{\rangle}
  \def\nn{\nabla} \def\pp{\partial} \def\E{\mathbb E}
\def\d{\text{\rm{d}}} \def\bb{\beta} \def\aa{\alpha} \def\D{\scr D}
  \def\si{\sigma} \def\ess{\text{\rm{ess}}}\def\s{{\bf s}}
\def\beg{\begin} \def\beq{\begin{equation}}  \def\F{\scr F}
\def\Ric{\mathcal Ric} \def\Hess{\text{\rm{Hess}}}
\def\e{\text{\rm{e}}} \def\ua{\underline a} \def\OO{\Omega}  \def\oo{\omega}
 \def\tt{\tilde}\def\[{\lfloor} \def\]{\rfloor}
\def\cut{\text{\rm{cut}}} \def\P{\mathbb P} \def\ifn{I_n(f^{\bigotimes n})}
\def\C{\scr C}      \def\aaa{\mathbf{r}}     \def\r{r}
\def\gap{\text{\rm{gap}}} \def\prr{\pi_{{\bf m},\varrho}}  \def\r{\mathbf r}
\def\Z{\mathbb Z} \def\vrr{\varrho} \def\ll{\lambda}
\def\L{\scr L}\def\Tt{\tt} \def\TT{\tt}\def\II{\mathbb I}
\def\i{{\rm in}}\def\Sect{{\rm Sect}}  \def\H{\mathbb H}
\def\M{\mathbb M}\def\Q{\mathbb Q} \def\texto{\text{o}} \def\LL{\Lambda}
\def\Rank{{\rm Rank}} \def\B{\scr B} \def\i{{\rm i}} \def\HR{\hat{\R}^d}
\def\to{\rightarrow}\def\l{\ell}\def\iint{\int}\def\gg{\gamma}
\def\EE{\scr E} \def\W{\mathbb W}
\def\A{\scr A} \def\Lip{{\rm Lip}}\def\S{\mathbb S}
\def\BB{\scr B}\def\Ent{{\rm Ent}} \def\i{{\rm i}}\def\itparallel{{\it\parallel}}
\def\g{{\mathbf g}}\def\Sect{{\mathcal Sec}}\def\T{\mathcal T}\def\BB{{\bf B}}
\def\f{\mathbf f} \def\g{\mathbf g}\def\BL{{\bf L}}  \def\BG{{\mathbb G}}
\def\Bd{{D^E}} \def\BdP{D^E_\phi} \def\Bdd{{\bf \dd}} \def\Bs{{\bf s}} \def\GA{\scr A}
\def\Bg{{\bf g}}  \def\Bdd{\psi_B} \def\supp{{\rm supp}}\def\div{{\rm div}}
\def\ddiv{{\rm div}}\def\osc{{\bf osc}}\def\1{{\bf 1}}\def\BD{\mathbb D}\def\GG{\Gamma}
\def\H{{\bf H}}
\maketitle

\begin{abstract} In recent years, remarkable progress has been made for Distribution dependent stochastic equations (DDSDEs) with singular interactions, existing results include  well-posedness, propagation of chaos, entropy-cost inequality and ergodicity. As a continuation to the existing study, 
in this paper we establish  Bismut type formulas   for the intrinsic derivative of  DDSDEs with singular interactions, which extends the existing formula established for the case with Lions differentiable drifts.     \end{abstract} \noindent


 \vskip 2cm

 \section{Introduction}

 Since 1996 when McKean  \cite{McKean} proposed a class of stochastic differential equations (SDEs) with   distribution dependent drifts to characterize nonlinear Fokker-Planck equations, distribution dependent SDEs  (also called McKean-Vlasov SDEs) have been intensively investigated.
In recent years, some remarkable progress has been made for SDEs with singular interactions, where the drift contains convolutions of distribution with a singular kernel, which describes singular interaction
in the corresponding mean-field particle system, see   \cite{CJM25,HRZ,HZ, HRW25,HRW26,IS,JW} and reference therein for the study on well-posedness, propagation of chaos, regularity estimates, and ergodicity.

On the other hand, as a powerful tool characterizing the regularity of stochastic systems, the derivative formula initiated by Bismut \cite{B} for the heat semigroup on manifolds
has been extensively studied and applied to various models. In recent years, Bismut type derivative formulas  with respect to the initial distribution have been established for distribution dependent SDEs,
see \cite{BRW, HSW21, HW*, RW19, W23} for the intrinsic/Lions derivative formulas, and see \cite{Ren} for the extrinsic derivative formula.
However, the existing Bismut formula for the intrinsic derivative  requires that the drift is either  Lions differentiable in the distribution variable \cite{W23},
or has half-Dini continuous extrinsic derivative \cite{HW*}.  These conditions excludes McKen-Vlasov SDEs with singular interactions, where the drift contains a term like
\beq\label{SG} B_t(x,\mu)=(h_t*\mu)(x):= \int_{\R^d} h_t(x-y)\mu(\d y),\ \ \ t\ge 0,\ x\in \R^d,\ \mu\in \hat{\scr P},\end{equation}
where
$$h: [0,\infty)\times \R^d \to \R^d$$ is measurable and singular in the spatial variable $x\in \R^d$, and $\hat{\scr P}$ the set  of probability measures on $\R^d$ such that the convolution
with $h_t$ is well-defined. Typical examples of $h_t$ include vector fields satisfying
$$|h_t(x)|\le \ff{c t^\tau }{|x|^\bb},\ \ \ x\ne 0,\ t\ge 0,$$
for some constants $\tau\ge 0,\bb\ge 0.$ It is the case when
$$h_t(x)= \ff {ct^\tau x}{|x|^{\bb+1}},\ \ \ x\in\R^d\setminus \{0\}, \ t\ge 0,$$ for which  $b_t(x,\mu)$ is not intrinsically  differentiable in $\mu$, and the extrinsic derivative is singular even if exists,
so that existing Bismut formula for the intrinsic derivative of distribution dependent SDEs do not work.

 In this paper,  we intend to establish the Bismut type derivative formula for the   following McKean-VlSDE on $\R^d$ with  singular interaction given by \eqref{SG}:
 \beq\label{E00} \d X_t= \big(b_t(X_t)+B_t(X_t,\L_{X_t})\big)\d t+\si_t(X_t)\d W_t,\ \ t\in [0,T],\end{equation}
 where $T>0$ is a fixed time, $W_t$ is the $d$-dimensional Brownian motion on a probability base $(\OO,\{\F_t\}_{t\in [0,T]},\F,\P)$, $\L_{X_t}$ is the distribution of $X_t$,
 and  $$b: [0,T]\times \R^d \to \R^d,\ \ B: [0,T]\times \R^d\times \hat{\scr P}\to \R^d,\ \ \si: [0,T]\times \R^d\to \R^d\otimes \R^d$$
 are measurable, where $\hat{\scr P}$ is equipped with the weak topology. When different probability space are concerned, we use $\L_{X_t|\P}$ to denote the distribution of $X_t$ under probability $\P$, and
 let $\L_{X_{[0,T]}|\P}$ be the distribution of the process $(X_t)_{t\in [0,T]}$ under $\P$.

 \beg{defn}\emph{ A continuous adapted process $(X_t)_{t\in [0,T]}$ on $\R^d$ is called a solution of \eqref{E00}, if $\P$-a.s.
 $$X_t= X_0+ \int_0^t  \big(b_s(X_t)+B_s(X_s,\L_{X_s})\big)\d s+\int_0^t \si_s(X_s)\d W_s,\ \ \ t\in [0,T],$$
 which means that the (stochastic) integrals in the right hand side are well-defined such that the equation holds.}

 \emph{A  couple $(X_t,W_t)_{t\in [0,T]}$ is called a weak solution of \eqref{E00}, if there exists a probability base such that $W_t$ is a $d$-dimensional Brownian motion and
 $X_t$ solves \eqref{E00}.}

 \emph{We call \eqref{E00} has a unique weak solution with initial distribution $\mu\in \scr P$, if the weak solution exists, and for any two
 weak solutions $(X_t^i,W_t^i)_{t\in [0,T]}$ with respect to probabilities $\P^i, i=1,2$ we have $\L_{X_{[0,T]}^1|\P^2}= \L_{X_{[0,T]}^2|\P^2}$.}
  \end{defn}

 Let $\scr P$ be  the set of all probability measures on $\R^d$ equipped with the weak topology. For any  $p\in [1,\infty)$, let
 $$\scr P_p:=\big\{\mu\in \scr P: \ \mu(|\cdot|^p)<\infty\big\}.$$
   If for any $\F_0$-measurable initial value $X_0$ with $\mu=\L_{X_0}\in \scr P_p$,
 the SDE \eqref{E00} has a unique solution, we denote
 $$P_t^*\mu=\L_{X_t},\ \ \ t\in [0,T].$$
 We study the intrinsic derivative of
 $$\mu\mapsto P_tf(\mu):= \int_{\R^d} f\d(P_t^*\mu)$$ for $t\in (0,T]$ and $f\in \B_b(\R^d),$
 where $ \B_b(\R^d)$ is the set of all bounded measurable functions on $\R^d$.

The intrinsic derivative for measures was introduced in \cite{AKR} to  construct diffusion processes on the configuration space  over   Riemannian manifolds,
and had been extended   in \cite{BRW}  for functions  on the $L^p$-Wasserstein space over a Banach space.

   For any $\mu\in \scr P_p$, the tangent space at $\mu$ is
 $$T_{\mu,p} := L^k(\R^d\to\R^d;\mu).$$
 Then $\mu\circ(id+\phi)^{-1} \in \scr P_p$ for $\phi\in T_{\mu,p}$, where $id$ is the identity map on $\R^d$.

\begin{defn} Let  $p\in (1,\infty)$ and $f$  be a continuous function  on $\scr P_p$.
It  is called intrinsically differentiable at a point $\mu\in\scr P_p$, if
 $$T_{\mu,p}\ni\phi\mapsto D_\phi f(\mu):= \lim_{\vv\downarrow 0} \ff{f(\mu\circ(id+\vv\phi)^{-1})-f(\mu)}{\vv}\in\R
$$ is a well  defined bounded linear functional. In this case,  the intrinsic derivative is the unique element
$$Df(\mu)\in T_{\mu,p}^*:=L^{\ff p{p-1}}(\R^d\to\R^d;\mu) $$ such that
$$\int_{\R^d}\<Df(\mu)(x), \phi(x)\>  \mu(\d x) = D_\phi f(\mu),\ \ \phi\in T_{\mu,p}.$$
 \end{defn}

 In Section 2, we state the main result of the paper (Theorem 2.1), which provides a Bismut-type formula for the intrinsic derivative of McKean–Vlasov SDEs with singular interactions under assumption (H).  In Section 3, we prove the first part of Theorem 2.1 (existence and moment estimate) by establishing well-posedness of the linearized derivative SDE (Lemma 3.1) and showing convergence of the solution process (Lemma 3.2). The estimates and techniques developed in Section 3 are then used in Section 4 to prove the second part of Theorem 2.1, namely the intrinsic Bismut formula (2.10) and its application to the semigroup (4.2).

 \section{Main result}

To characterize  the singularity of coefficients $b$ and $\si$  in time-space variables, we recall some functional spaces introduced in \cite{XXZZ}.
For any $p,q\ge1$, let $\tt L_q^p$ denote the set of measurable functions $f$ on $[0,T]\times\R^d$ such that
$$\|f\|_{\tt L_q^p}:= \sup_{z\in \R^d}\bigg( \int_{0}^{T} \|1_{B(z,1)}f_t\|_{L^p(\R^d)}^q\d t\bigg)^{\ff 1 q}<\infty,$$
where $B(z,1):= \{x\in\R^d: |x-z|\le 1\}$, and $\|\cdot\|_{L^p(\R^d)}$ is the $L^p$-norm with respect to the Lebesgue measure on $\R^d$.
We will  take $(p,q)$ from the class
$$\scr K:=\Big\{(p,q): p,q\in  (2,\infty),\  \ff d p+\ff 2 q<1\Big\}.$$ \index{$\scr K:=\{(p,q): p,q\in (1,\infty), \ff d p+\ff 2 q<1\}$}
According to \cite{HRW25}, for any $p\in [1,\infty]$,
$$\scr P_{p*}:=\bigg\{\mu\in \scr P:\ \|\mu\|_{p*}:= \sup_{f\in \B_b(\R^d), \|f\|_{\tt L^p}\le 1} |\mu(f)|<\infty\bigg\},$$
  is a complete metric space under the  $p*$-distance   $\|\mu-\nu\|_{p*}.$ It is clear that when $p=\infty$ we have
  $\scr P_{\infty*}=\scr P$ and $\|\mu-\nu\|_{\infty*}=\|\mu-\nu\|_{var}$ is the variation distance.
We make the following assumption on the coefficients of \eqref{E00}.

 \beg{enumerate}\item[$(H)$]    Let  $a_t(x):= (\si_t\si_t^*)(x)$ and $B_t(x,\mu)$ be given in $\eqref{SG}$,  where $\si^*$ is the transposition of $\si$.
  \item[$(1)$]  $a_t(x)$ is invertible with $\|a\|_\infty+\|a^{-1}\|_\infty<\infty$,    and    
  $$\|\nn \si \|\le \sum_{i=1}^l f_i $$    for some  $0\le f_i\in \tt L_{q_i}^{p_i}$ and $(p_i,q_i)\in \scr K,\ 1\le i\le l.$
 \item[$(2)$]    $b_t(x)$ satisfies
$$ \sup_{t\in [0,T],\ x\ne y} \bigg(|b_t(0)|+ \ff{|b_t(x)-b_t(y)|}{|x-y|}\bigg)<\infty.$$
  \item[$(3)$] For any $t\in [0,T]$, $\nn h_t\in \tt L^{k'}$ and is a.e. continuous. There exist constants $\kk\ge 0, K>0, k\in (d,\infty]$ and $k'\in (1,\infty]$ such that
$$\|h_t\|_{\tt L^{k}} \le K,\ \ \
  \big\|\nn h_t\big\|_{\tt L^{k'}} \le K t^\kk,\ \ t\in (0,T].$$
    \end{enumerate}

Since $(H)$ implies $(A_1)$ and $(A_2)$ in \cite{HRW25} for $p=\infty$, according to \cite[Theorem 2.1(2)-(3)]{HRW25}, for any initial value (respectively, initial distribution), $(H)$ implies that \eqref{E00} has a unique solution satisfying
\beq\label{*0} \sup_{t\in (0,T]} t^{\ff d{2k}} \|\L_{X_t}\|_{k*}<\infty,\ \ \ \E^{\F_0} \bigg[\sup_{t\in [0,T]} |X_t|^n\bigg]\le C_n(1+|X_0|^n),\ \ n\in [1,\infty)\end{equation}
for some $C_n\in (0,\infty),$ where $\E^{\F_0}$ is the conditional expectation given $\F_0$.

For any $\mu\in \scr P$, let $X_t^\mu$ be the unique solution with initial distribution $\mu$, and denote
$$P_t^*\mu:=\L_{X_t^\mu},\ \ \ \ t\in [0,T], \ \mu\in \scr P.$$
Since $\|\mu\|_{p*}=1$ for $p=\infty$, \cite[(2.2)]{HRW25} implies
\beq\label{MM} \sup_{t\in (0,T],\mu\in \scr P} t^{\ff d{2k}} \|P_t^*\mu\|_{k*} <\infty.\end{equation}
  We intend to establish the Bismut type formula for the intrinsic derivative of
$$\mu\mapsto P_tf(\mu):=\int_{\R^d} f\d (P_t^*\mu),\ \ \ \ t\in (0,T],\ f\in \B_b(\R^d).$$

To this end, given $\mu\in \scr P$ we   recall the Bismut formula for  the decoupled SDE
\beq\label{DC} \d X_t^{\mu,x}= \big\{b_t(X_t^{\mu,x})+ B_t(X_t^{\mu,x},P_t^*\mu)\big\}\d t+ \si_t(X_t^{\mu,x})\d W_t,\ \ \ (t,x)\in [0,T]\times\R^d,\ X_{0}^{\mu,x}=x.\end{equation}
By $(H)$(3) and \eqref{MM}, $b_t^\mu(x):= B_t(x, P_t^*\mu)$ satisfies
$$|b_t^\mu(x)|\le K\|P_t^*\mu\|_{k*}\le c t^{-\ff d{2k}},\ \ \ t\in (0,T],$$
which together with $k>d$ due to $(H)$(3) yields
\beq\label{BM}\sup_{\mu\in \scr P}  \|b^\mu\|_{\tt L_{q_0}^{p_0}}<\infty\ \text{ for \ some\ } (p_0,q_0)\in \scr K. \end{equation}
Combining this with $(H)$(1)-(2), we see that conditions $(A^{1.1}), (A^{1.2})$ and $(A^{1.3})$ in \cite{WR25} hold, so that
by \cite[Theorem 1.3.1 and Theorem 1.4.2]{WR25},  see also \cite[Theorem 2.1]{W23},
$$\nn_vX_t^{\mu,x}:=\lim_{\vv\downarrow 0} \ff{X_t^{\mu,x+\vv v}-X_t^{\mu,x}}\vv,\ \ \ t\in [0,T],\ x,v\in\R^d$$
exists in $L^n(\OO\mapsto\R^d,\P)$ for any $n\in [1,\infty)$,  and
there exists $C: [1,\infty)\to (0,\infty)$ such that
\beq\label{DS}  \E\bigg[\sup_{t\in [0,T]} |X_t^{\mu,x}|^n\bigg]\le C(n) (1+|x|^n),\ \ \
\E\bigg[\sup_{t\in [0,T]} |\nn_v X_t^{\mu,x}|^n\bigg]\le C(n) |v|^n\end{equation}  holds for all $\mu\in \scr P,\ x,v\in \R^d,$ and $ n\in [1,\infty).$
Moreover,
 $$P_t^\mu f(x):= \E\big[f(X_t^{\mu,x})\big],\ \ \ t\in [0,T],\ x\in\R^d,\ f\in \B_b(\R^d)$$
 satisfies the Bismut formula: for any $t\in (0,T] $ and $\bb\in C^1([0,1])$ with $\bb_0=0$ and $\bb_t=0$,
 \beq\label{B0} \nn_v P_t^\mu f(x)= \E\bigg[f(X_t^{\mu,x})\int_0^t \bb_s'\big\<(\si_sa_s^{-1})(X_s^{\mu,x})\nn_v X_s^{\mu,x}, \d W_s\big\>\bigg],\ \ \ f\in \B_b(\R^d),\ x,v\in\R^d.\end{equation}
 Finally, by $(H)$(1)-(2) and \eqref{BM}, the conditions $(C)$ and $(A_2)$ in \cite{HRW25} hold,
 so  by \cite[Proposition 5.4]{HRW25}, for any $1<p_1\le p_2\le\infty$, there exists a constant $c(p_1,p_2)\in (0,\infty)$ such that
 \beq\label{GRD}  \|\nn^i P_t^\mu\|_{\tt L^{p_1}\to \tt L^{p_2}}\le c(p_1,p_2) t^{-\ff{i} 2-\ff{d(p_2-p_2)}{2p_1p_2}},\ \ \ t\in (0,T],\ \mu\in \scr P,\ i=0,1,\end{equation}
 where $\nn^0$ is the identity map and $\nn^1=\nn$ is the gradient.

 For any $\F_0$-measurable random variable  $\eta$ on $\R^d$, consider the distribution dependent SDE
 $$\d X_t^{\eta,\mu} = \big\{b_t(X_t^{\eta,\mu})+B_t(X_t^{\eta,\mu},\L_{X_t^{\eta,\mu}})\big\}\d t+ \si_t(X_t^{\eta,\mu})\d W_t,\ \ t\in [0,T],\ X_0^{\eta,\mu}=X_0^\mu+\eta,$$
 and define
 $$\nn_\eta X_t^\mu:= \lim_{\vv\downarrow 0} \ff{X_t^{\vv\eta,\mu}-X_t^\mu}\vv$$
 if the limit exists in $L^p(\OO\to\R^d,\P)$ for some $p\in [1,\infty).$

We are now ready to state the  Bismut type formula for the intrinsic derivative of $P_t f(\mu)$.

\beg{thm}\label{T1} Assume $(H)$ with $2\kk-\ff d{k'}>-1$, and let $p\in [\ff{2k'(\kk+1)}{2k'(\kk+1)-d},\infty)\cap (\ff {k'}{k'-1},\infty)$.  Then the following assertions hold.
\beg{enumerate}  \item[$(1)$] There exists a constant $c\in(0,\infty)$  such that
\beq\label{Y1} \bigg(\E\Big[\big\|\nn h_t(z-X_t^\mu)\big\|^{\ff{p}{p-1}}\Big]\bigg)^{\ff{p-1}p}\le c t^{\kk-\ff{d p}{2k'(p-1)}},\ \ \ t\in (0,T],\ z\in\R^d,\ \mu\in \scr P.\end{equation}
Moreover, for any $\eta\in L^p(\OO\to\R^d,\F_0,\P)$,
$\nn_\eta X_t^\mu$ exists in $L^p(\OO\to\R^d,\P)$  and there exists a constant $c\in (0,\infty)$ independent of $\mu$ and $\eta$, such that
\beq\label{YB} \sup_{t\in [0,T]}\E\Big[  \big|\nn_\eta X_t^\mu\big|^p\Big]\le c\E[|\eta|^p].\end{equation}
\item[$(2)$]   Denote $\zeta=\si(\si\si^*)^{-1}$. For any $\phi\in T_{\mu,p}, f\in \B_b(\R^d)$, $t\in (0,T]$ and  and $\bb\in C^1([0,1])$ with $\bb_0=0$ and $\bb_t=0$,
$D_\phi P_t f(\mu)$ exists and
\beq\label{BSMI} \beg{split} & D_\phi P_tf(\mu)= \int_{\R^d} \E\bigg[f(X_t^{\mu,x}) \int_0^t \bb_s'\big\<\zeta_s(X_s^{\mu,x} ) \nn_{\phi (x)}X_s^{\mu,x},  \d W_s\big\> \bigg]\mu(\d x)\\
&+  \E\bigg[f(X_t^{\mu}) \int_0^t  \Big\<\zeta_s(X_s^\mu)  \E\big[\<\nn h_s  (z -X_s^\mu), \nn_{\phi(X_0^\mu)} X_s^\mu\>\big]|_{z=X_s^\mu},  \d W_s\Big\> \bigg].
\end{split} \end{equation}
Consequently,   there exists a constant $c>0$ such that
$$\|D P_tf(\mu)\|_{L^{\ff p{p-1}}(\mu)} \le \ff {c}{\ss t} \big\|f(X_t^\mu)\big\|_{L^{\ff p{p-1}}(\P)},  \
  \ t\in (0,T], f\in \B_b(\R^d),\mu\in \scr P_p.$$
\end{enumerate}
\end{thm}

By Zvonkin's transform, one may allow the drift also contains local integrable terms in time-spatial, but the integrable indexes of these terms will lead restriction
to the index $p$. Since our main concern is about the interaction kernel $h_t$,  in this paper  we do not consider this situation.
Below we present an example to show that our result applies to irregular kernel $h_t$ such that the convolution is not differentiable in distribution.

\paragraph{Example 2.1} Let $d\ge 2$. Consider the SDE \eqref{E00} with singular interaction where $B_t(x,\mu)$ is in \eqref{SG} with $h_t$ satisfying
$$|h_t(z)|\le  C t^\kk |z|^{-\bb},\ \ \ \ t\in [0,T],\ 0\ne z\in\R^d$$
for some constants $\bb\in (0,1)$ and $\kk\in (\ff\bb 2,\infty)$. Then $(H)$(3) holds for $k\in (d, \ff d \bb)$ and $k'\in (1\lor \ff d{2\kk+1}, \ff d{\bb+1}).$
Therefore, if $\si$ and $b$ satisfy $(H)$(1)-(2), and $\kk>\ff \bb 2$,    assertions in  Theorem \ref{T1} hold for any $p\in [\ff{2k'(\kk+1)}{2k'(\kk+1)-d},\infty)\cap(\ff{k'}{k'-1},\infty)$.

\section{Proof of Theorem \ref{T1}(1)}
By $(H)$(3) and \eqref{GRD} for $i=0$, we  find a constant $c\in (0,\infty)$ independent of $z$ and $\mu$ such that
\beg{align*} &\big(\E\|\nn h_t\|(z-X_t^\mu)\|^{\ff p{p-1}}\big)^{\ff{p-1}p} = \Big[\mu\Big(P_t^\mu \|\nn h_t\|(z-\cdot)\|^{\ff p{p-1}}\Big)\Big]^{\ff{p-1}p}\\
&\le \|P_t^\mu\|_{\tt L^{\ff{k'(p-1)} p}\to\tt L^\infty} \big\|\nn h_t(z-\cdot)\big\|_{\tt L^{k'}} \le c t^{\kk- \ff{dp}{2k'(p-1)}},\ \ \ t\in (0,T].\end{align*}
So, \eqref{Y1} holds.  It remains to prove the  existence of $\nn_\eta X_t^\mu$ and to verify the estimate \eqref{YB}.
This follows from the following two lemma.

In the following, let $\mu\in \scr P$ and $\eta\in L^p(\OO\to\R^d,\F_0,\P)$ be fixed,
we simply denote
$$X_t=X_t^\mu,\ \ \  \mu_t= P_t^*\mu,\ \ \ t\in [0,T].$$  By formal calculations,  if $\nn_\eta X_t^\mu$ exists then it satisfies the following SDE:
\beq\label{VT}\beg{split}&\d v_t= \Big\{\nn_{v_t} b_t(X_t) + \nn_{v_t} (h_t*\mu_t) (X_t) +\E\big[\nn_{v_t}h_t(z-X_t)\big]\big|_{z=X_t}\Big\}\d t\\
&\qquad \qquad \qquad + \nn_{v_t}\si_t(X_t)\d W_t,
\ \ t\in [0,T],\ v_0=\eta.\end{split} \end{equation}
So, we first solve this SDE and estimate the moments of the solution.

\beg{lem}\label{L1} Assume $(H)$ with $2\kk-\ff d{k'}>-1$ and let $p\in [\ff{2k'(\kk+1)}{2k'(\kk+1)-d},\infty)$.  Then for any $\mu\in \scr P$ and $\eta\in L^p(\OO\to\R^d,\F_0,\P)$,  the  SDE $\eqref{VT}$ has a unique solution
satisfying
\beq\label{NB} \rr_T(v):=   \bigg(\E\Big[\sup_{t\in [0,T]}   |v_t|^p\Big]\bigg)^{\ff 1 p}<\infty.\end{equation}
 Moreover, for any $n\ge 1$ there exists a constant $c\in (0,\infty)$ independent of $\mu$ and $\eta$, such that
\beq\label{VT'} \sup_{t\in [0,T]}\E\big[ |v_t|^n\big|\F_0\big]\le c|\eta|^n.\end{equation}\end{lem}

\beg{proof}  We first verify \eqref{VT'}, then prove the well-posedness for bounded $\|\nn \si\|$, and finally make extension by mollifier approximation.

(a) Proof of  \eqref{VT'}.  Let $v=(v_t)_{t\in [0,T]}$  solve \eqref{VT} satisfying \eqref{NB}. Since $p\in [\ff{2k'(\kk+1)}{2k'(\kk+1)-d},\infty)$, we have
\beq\label{KL} \kk-\ff{dp}{2k'(p-1)}>-1.\end{equation}
By \eqref{Y1} and H\"older's inequality, we obtain
\beq\label{X1} \beg{split}&\big| \nn_{v_s} (h_s*\mu_s) (z)\big|+ \Big|\E\big[\nn_{v_s}h_s(z-X_s)\big] \Big|\\
&\le \Big( |v_s| +  \big(\E|v_s|^p\big)^{\ff 1 p} \Big) \Big(\E\big[|\nn h_s(z-X_s)|^{\ff p{p-1}}\big]\Big)^{\ff{p-1}p}\\
&\le c \Big( |v_s| +  \big(\E|v_s|^p\big)^{\ff 1 p} \Big) s^{\kk-\ff {dp}{2k'(p-1)}},\ \ \ s\in (0,T],\ z\in\R^d.\end{split}\end{equation}
Combining this with $(H)$(1)-(2) and  It\^o's formula, for any $n\ge p$, we find a constant $c_1(n)\in (0,\infty)$ such that
$$\d |v_t|^{2n} \le c_1(n) \Big[|v_t|^{2n} \big(t^{\kk-\ff {dp}{2k'(p-1)}}+\|\nn \si_t(X_t)\|^2\big)+ t^{\kk-\ff {dp}{2k'(p-1)}}\big(\E|v_t|^p\big)^{\ff {2n} p}\Big]\d t + \d M_t$$
for some martingale $M_t$. By the stochastic Gronwall inequality for the conditional expectation $\E^{\F_0}:=\E(\cdot|\F_0)$, see for instance \cite[Lemma 1.3.3]{WR25}
for $q=\ff 1 2$ and $p=\ff 3 4$,   we find a constant $c_2\in (0,\infty)$ such that
$$\E^{\F_0}\bigg[\sup_{s\in [0,t]} |v_s|^n\bigg]  \le c_2 \E^{\F_0}\Big[\e^{c_2 \int_0^t \|\nn \si_s(X_s)\|^2 \d s}\Big] \bigg(|\eta|^{2n} +
\int_0^t s^{\kk-\ff {dp}{2k'(p-1)}}\big(\E|v_s|^p\big)^{\ff {2n} p}\d s\bigg)^{\ff 1 2}.$$
Combining this with the Krylov and Khasminskii estimates, see for instance \cite[Theorem 1.2.3(2) and Theorem 1.2.4]{WR25},
we find a constant $c_3\in (0,\infty)$ such that
\beq\label{POO}  \E^{\F_0}\bigg[\sup_{s\in [0,t]} |v_s|^n\bigg]  \le c_3|\eta|^{n} +
c_3\bigg(\int_0^t s^{\kk-\ff {dp}{2k'(p-1)}}\big(\E|v_s|^p\big)^{\ff {2n} p}\d s\bigg)^{\ff 1 2}.\end{equation}
Taking $n=p$ this implies
$$ \E\big[|v_t|^p\big]\le \E\bigg[\E^{\F_0} \bigg(\sup_{s\in [0,t]} |v_s|^p\bigg)\bigg]\le c_3\E\big[|\eta|^p\big]+ c_3
\bigg(\int_0^t s^{\kk-\ff {dp}{2k'(p-1)}}\big(\E|v_s|^p\big)^{2}\d s\bigg)^{\ff 1 2}.$$
By \eqref{KL},  this together with \eqref{NB}  implies
$$\E\big[|v_s|^p\big]\le c(p)(1+\E[|\eta|^p]),\ \ \ s\in [0,T]$$
for some constant $c(p)\in [0,\infty].$ Substituting into \eqref{POO} we prove  \eqref{VT'}.

(b) The well-posedness of \eqref{VT} for bounded $\|\nn\si\|$.
Let  $\scr V_T$ be the space of all adapted  continuous process $v=(v_t)_{t\in [0,T]}$ on $\R^d$ satisfying \eqref{NB}.
For any $\ll\in [0,\infty)$, $\scr V_T$  is a complete space under the metric
$$\rr_{T,\ll}(v, w):=    \sup_{t\in [0,T]}\e^{-\ll t}\bigg(\E\Big[\sup_{s\in [0,t]} |v_s-w_s|^p\Big]\bigg)^{\ff 1p},\ \ \ v, w \in \scr V_T.$$
 For any $v\in \scr V_T,$ define $\psi(v)=(\psi_t(v))_{t\in [0,T]}$ by
\beq\label{GG1} \beg{split} \psi_t(v):= \eta +&\int_0^t \Big\{\nn_{v_s} b_s(X_s) + \nn_{v_s} (h_s*\mu_s) (X_s) +\E\big[\nn_{v_s}h_s(z-X_s)\big]\big|_{z=X_s}\Big\}\d s\\
 & + \int_0^t\nn_{v_s}\si_s(X_s)\d W_s,\ \ \ t\in [0,T].\end{split}\end{equation}
By $(H)$(2) and \eqref{NB} for $w$ in place of ,  there exists a constant $c_1\in (0,\infty) $ such that
 \beq\label{GG2}  \E \int_0^T  |\nn_{v_s} b_s(X_s)|^p\d s  \le c_1 \int_0^T   \E\big[|v_s|^p\big]\d s \le c_1 \rr_T(v)^p<\infty.\end{equation}
 Noting that for any $z\in\R^d$,
 \beg{align*} &\big|\nn_{v_s} (h_s*\mu_s)(z)\big|\le |v_s| \E\big[|\nn h_s(z-X_s)|\big]\le |v_s| \Big(\E\big[|\nn h_s(z-X_s)|^{\ff p{p-1}}\big]\Big)^{\ff{p-1}p},\\
 &\big|\E\big[\nn_{v_s} h_s(z-X_s)]\big|\le \Big(\E\big[|v_s|^p\big]\Big)^{\ff 1 p}  \Big(\E\big[|\nn h_s(z-X_s)|^{\ff p{p-1}}\big]\Big)^{\ff{p-1}p},\end{align*}
 by \eqref{Y1}  and \eqref{KL},  we   find a constant $c_2\in (0,\infty)$ such that
  \beq\label{GG3}\beg{split} &\E\bigg(\int_0^T \Big(\big| \nn_{v_s} (h_s*\mu_s) (X_s)\big|+ \Big|\E\big[\nn_{v_s}h_s(z-X_s)\big]\big|_{z=X_s}\Big|\Big)\d s\bigg)^p\\
  &\le \sup_{z\in\R^d}  \E\bigg(\int_0^T \Big( |v_s| +  \big(\E|v_s|^p\big)^{\ff 1 p} \Big) \Big(\E\big[|\nn h_s(z-X_s)|^{\ff p{p-1}}\big]\Big)^{\ff{p-1}p}\d s \bigg)^p\\
  &\le  c_2\E\bigg[\sup_{t\in [0,T]} |v_s|^p\bigg]  \bigg(\int_0^T  s^{\kk-\ff{dp}{2k'(p-1)}}\d s\bigg)^p<\infty.\end{split}\end{equation}
Moreover, by the Burkholder-Davis-Gundy inequality, $\|\nn \si\|_\infty<\infty$ as assumed in the moment,  and   \eqref{NB},  we find a constant  $c_3 \in (0,\infty)$ such that
\beg{align*} \E\bigg[\sup_{t\in [0,T]}\bigg|\int_0^t\nn_{v_s}\si_s(X_s)\d W_s\bigg|^p\bigg] \le c_3 \E\bigg(\int_0^T |v_s|^2  \d s\bigg)^{\ff p2} <\infty.\end{align*}
Combining this  with \eqref{GG1}-\eqref{GG3},  we obtain $\psi(v)\in \scr V_T$, so that  $\psi: \scr V_T\to \scr V_T.$

Similarly, by   $(H)$, \eqref{Y1}    and Jensen's inequality,   we find constants $c_4,c_5,c_6\in (0,\infty)$ such that for any  $v,w\in \scr V_T,$
 \beg{align*}&\E\bigg[ \sup_{s\in [0,t]} |\psi_s(v)-\psi_s(w)|^p\bigg]\\
 & \le c_4  \E \bigg(\int_0^t  \Big[|v_s-w_s| + \big(\E|v_s-w_s|^p \big)^{\ff 1 p} \Big]s^{\kk-\ff{dp}{2k'(p-1)}} \d s\bigg)^p
+ c_5 \E\bigg(\int_0^t  |v_s- w_s|^2\d s\bigg)^{\ff p 2}\\
&\le c_6 \int_0^t \E\big[|v_s-w_s|^p\big]  s^{\kk-\ff{dp}{2k'(p-1)}} \d s+ c_5 \E\bigg(\int_0^t  |v_s- w_s|^2\d s\bigg)^{\ff p 2}.  \end{align*}
Combining this with \eqref{KL}, we conclude that for  $\ll\in (0,\infty)$   large enough,
\beq\label{CT} \beg{split}  &\rr_{T,\ll}\big(\psi(v), \psi(w)\big)^p= \sup_{t\in [0,T]} \e^{-\ll t} \E\bigg[ \sup_{s\in [0,t]} |\psi_s(v)-\psi_s(w)|^p\bigg]\\
&\le  \rr_{T,\ll}(v,w)^p \sup_{t\in [0,T]} \bigg[\int_0^t \e^{-\ll p(t-s) } s^{\kk-\ff{dp}{2k'(p-1)}} \d s+  \bigg(\int_0^t  \e^{-2\ll  (t-s)} \d s\bigg)^{\ff p 2}\bigg]\\
&\le \ff 1 2   \rr_{T,\ll}(w,\tt w)^p,\ \ \  v, w\in  \scr V_T.\end{split}\end{equation}
By the Banach fixed point theorem, this implies that $\psi$ has a unique fixed point
 $v\in\scr V_T$, which is the unique solution to \eqref{VT} satisfying \eqref{NB}.

(c) Well-posedness of \eqref{VT} under $(H)$. We will make   approximations of $\nn \si_t$ as follows. 
 For any $m\in\mathbb N$, define the linear operator 
  $Q_m(t): \R^d\to \R^d\otimes \R^d$  satisfying 
 $$\big( Q_m(t)  e_k\big)_{ij} = \max\Big(-m, m\land \big(\pp_k (\nn \Theta_t )\si_t\big)_{ij}(X_t)\Big),\ \ \ 1\le i,j,k\le d, $$
 where   $\{e_k\}_{1\le k\le d}$ is the standard ONB of $\R^d$. 
 Then $\lim_{m\to\infty}Q_m(t)=  \big(\nn  \si_t \big)(X_t)$ and by $(H)$, there exists a constant $c>0$ such that 
\beq\label{PLO} \sup_{m\ge 1}  \|Q_m(t)\|^2\le  \le c+ c\sum_{i=0}^l  f_i(t,X_t)^2=:f(t,X_t)^2.\end{equation}  
So, by Krylov's estimate and dominated convergence theorem, 
 \beq\label{PLLY}\lim_{m\to\infty} \E\bigg(\int_0^T \Big\|Q_m(t)- \nn \big\{\pp_k  \si_t \big\}(X_t)\Big\|^2\d s\bigg)^p=0,\ \ \ p\in [1,\infty).\end{equation} 
 By (a) and (b),   the SDE \eqref{VT} with $Q_n(t)  $ in place of $\nn \si_t(X_t)$ has a unique solution $v_t^n$ satisfying \eqref{VT'}.
For any $m,m'\ge 1$, we have
\beg{align*}&\d \big(v_t^m-v_t^{m'}\big)= \Big(\nn_{v_t^m-v_t^{m'}} b_t(X_t^\mu) + \nn_{v_t^m-v_t^{m'}} (h_t*\mu_t) (X_t) +\E\big[\nn_{v_t^m-v_t^{m'}}h_t(z-X_t)\big]\big|_{z=X_t}\Big)\d t\\
&\qquad  + \Big[\nn_{v_t^m-v_t^{m'}}Q_m(t) + \nn_{v_t^{m'}} \big(Q_m(t) -Q_{m'}(t)\big)  \Big]\d W_t,
\ \ t\in [0,T],\ v_0^m-v_0^{n}=0.\end{align*}
So, by $(H)$(1)-(2),  \eqref{X1} for $v_t$ in place of $v_t^m-v_t^{m'}$ and \eqref{PLO}, for any $n\ge p$ we find a constant $C_1\in (0,\infty)$ such that
\beg{align*} &\d \big|v_t^m-v_t^{m'}\big|^{2n} \le C_1\Big( 1+ t^{\kk-\ff{dp}{2k'(p-1)}}+ f(t,X_t)^2 \|^2\Big)(X_t) \big|v_t^m-v_t^{m'}\big|^{2n}\d t\\
&\qquad + \Big[\big(\E \big|v_t^m-v_t^{m'}\big|^{p}\big)^{\ff {2n}p} t^{\kk-\ff{dp}{2k'(p-1)}}+ \|v_t^{n}\|^{2n} \|Q_m(t)-Q_{m'}(t)\|^{2} \Big]\d t+\d M_t,
\ \ t\in [0,T].\end{align*}
By the stochastic Gronwall inequality and Khasminskii inequality as in step (a),  we find   constants $C_2,C_3\in (0,\infty)$ such that
 \beg{align*} &\bigg(\E^{\F_0}\bigg[\sup_{s\in [0,t]} \big|v_s^m-v_s^{m'}\big|^{n}\bigg]\bigg)^2\\
 &\le C_2  \int_0^t  \Big(\big(\E \big|v_s^m-v_s^{m'}\big|^{p}\big)^{\ff {2n}p} s^{\kk-\ff{dp}{2k'(p-1)}}+ \E^{\F_0}\big[\|v_s^{m'}\|^{2n} \|Q_m(t)-Q_{m'}(t)\|^{2}\big] \Big)\d s.
   \end{align*}
 Moreover, by \eqref{VT'} for $v_t^{m'}$ in place of $v_t$, and by Krylov's estimate, see for instance \cite[Theorem 1.3.3(2)]{WR25}, we find a constant $C_3\in (0,\infty)$ and some
  positive constants  $\vv_{m,m'}\to 0$ as $m,m'\to\infty$ due to   \eqref{PLO} and \eqref{PLLY}, 
 such that
 \beg{align*} & \E^{\F_0} \int_0^t    \|v_s^{m'}\|^{2n} \|\nn(\si_s^m-\si_s^{m'})\|^{2} \d s\\
 &\le C_2   \bigg(\E^{\F_0}  \Big[\sup_{s\in [0,t]}\|v_s^{m'}\|^{4n}\bigg)^{\ff 1 2}
 \bigg(\E\bigg|\int_0^t \|Q_m(t)-Q_{m'}(t))\|^{2}\d s\bigg|^2\bigg)^{\ff 1 2}\\
 &\le C_3 |\eta|^{2n} \vv_{m,m'},\ \ \ t\in [0,T],\ m,m'\ge 1.
  \end{align*} By choosing  $n=p$ and taking expectation, we find a constant $C\in (0,\infty)$ such that
  \beg{align*}&\E\bigg[\sup_{s\in [0,t]} \big|v_s^m-v_s^{m'}\big|^{p}\bigg]\\
  &\le C  \bigg(\int_0^t  \Big(\big(\E \big|v_s^m-v_s^{m'}\big|^{p}\big)^{2} s^{\kk-\ff{dp}{2k'(p-1)}}\d s\bigg)^{\ff 1 2}
   + C  \vv_{m,m'},\ \ \ t\in [0,T],\ m,m'\ge 1.\end{align*}
By combining this with $\kk-\ff{dp}{2k'(p-1)}>-1$,  $v_t^m, v_t^{m'}\in \scr V_T$ so that $\rr_T(w^m-v^{m'})<\infty,$ and $\vv_m+\vv_{m'}\to 0$ as $m,m'\to\infty$,
we conclude that $\{v^m\}_{m\ge 1}$ is a Cauchy sequence in $\scr V_T$ and the limit $v=\lim_{m\to\infty}v^m$ is a solution of \eqref{VT}.

The uniqueness follows by applying It\^o's formula and stochastic Gronwall's inequality to  $|v_t-\tt v_t|^{2p}$ as above for two  solutions $v_t$ and $\tt v_t$. 

   \end{proof}

\beg{lem}\label{L2}   Assume $(H)$ with $2\kk-\ff d{k'}>-1$ and let $p\in [\ff{2k'(\kk+1)}{2k'(\kk+1)-d},\infty)$.  Let $v_t$ solve $\eqref{VT}$ satisfying $\eqref{NB},$  and let
   $v_t^\vv:= \ff{X_t^{\vv\eta,\mu}-X_t^\mu}\vv$ for $\vv>0,\  t\in [0,T].$  Then
$$\lim_{\vv\downarrow 0}\E\bigg[\sup_{t\in [0,T]} \big|v_t^\vv- v_t\big|^p\bigg]=0.$$
\end{lem}

\beg{proof} We complete the proof by three steps.

(a) We claim that for any $n\in [1,\infty)$ there exists a constant $K_n\in (0,\infty)$ such that
\beq\label{CE} \E^{\F_0} \bigg[\sup_{t\in [0,T]} |v_t^\vv|^n\bigg]\le K_n |\eta|^n,\ \ \ \vv\in (0,1].\end{equation}
To this end, denote
$$X_t^\vv:= X_t^{\vv\eta,\mu},\ \ \ \mu_t^\vv:= \L_{X_t^{\vv}} = P_t^* \L_{X_0^\mu+\vv \eta},$$ and let $(Y_t^\vv)_{t\in [0,T]}$ solve the SDE
$$ \d Y_t^\vv = \big\{b_t(Y_t^\vv)+ B_t(Y_t^\vv, \mu_t)\big\}\d t+\si_t(Y_t^\vv)\d W_t,\ \ \ Y_0^\vv= X_0^\vv= X_0^\mu+\vv\eta,\ \ \ t\in [0,T].$$
By $\|\nn b\|_\infty<\infty$ due to $(A)$, we find a martingale $M_t$ and a constant $D_n\in (0,\infty)$ such that
\beq\label{MFF} \beg{split} & \d |Y_t^\vv- X_t^\vv|^{2n} \le \d M_t +   2n(2n-1) |X_t^\vv-Y_t^\vv|^{2n-2} \|\si_t(Y_t^\vv)-\si_t(X_t^\vv)\|_{HS}^2 \d t \\
&  + D_n |X_t^\vv-Y_t^\vv|^{2n-1} \big(1+ |(h_t*\mu_t^\vv)(Y_t^\vv)- (h_t*\mu_t^\vv)(X_t^\vv)+ (h_t*(\mu_t-\mu_t^\vv)) (Y_t^\vv)\big| \big)\d t.
 \end{split}\end{equation}
 By \eqref{X1}, we obtain
\beq\label{MF0}\beg{split} & \big|(h_t*\mu_t^\vv)(Y_t^\vv)- (h_t*\mu_t^\vv)(X_t^\vv)\big|\le \int_0^1 \Big|\nn_{Y_t^\vv-X_t^\vv} h_t\big(rX_t^\vv+(1-r)Y_t^\vv\big)\Big| \d r \\
&\le c \Big(|X_t^\vv-Y_t^\vv|+ \big(\E[|X_t^\vv-Y_t^\vv|^p]\big)^{\ff 1 p}\Big) t^{\kk-\ff {dp}{2k'(p-1)}},\ \ \ t\in (0,T].\end{split}\end{equation}
Next, noting that
$$\W_1(\mu_0,\mu_0^\vv)\le \E[|\vv \eta|]\le \vv \big(\E[|\eta|^p]\big)^{\ff 1 p}<\infty,\ \ \vv\in (0,1],$$
by $(H)$(3) and \cite[Theorem 2.3(1)]{HRW25} for $p=\infty$ and $q=1$, we find a constant $c_0\in (0,\infty)$ such that
\beq\label{MF1} \big\|h_t*(\mu_t-\mu_t^\vv)\big\|_\infty \le \|h_t\|_{\tt L^k} \|\mu_t-\mu_t^\vv\|_{k*} \le c_0 \vv t^{-\ff 1 2-\ff d{2k}},\ \ \ t\in (0,T].\end{equation}
Moreover, by the maximal function estimate \cite[Lemma 2.1]{XXZZ}, we find a constant $c_1\in (0,\infty)$ such that
\beq\label{MF} \beg{split}& \|f(x)-f(y)\|\le c_1|x-y|\big(\|f\|_\infty +\scr Mf(x)+ \scr Mf(y)\big),\\
&\|\scr M f\|_{\tt L_q^p}\le c_{p,q} \|f\|_{\tt L_q^p},\ \ \ p,q\in (1,\infty),\ x,y\in\R^d,\end{split} \end{equation}
where $f$ is a nonnegative continuous function on $\R^d$, $c_{p,q}\in (0,\infty)$, and
$$\scr M f(x):=\sup_{r\in (0,1)} \ff 1 {|B(0,r)|}\int_{B(0,r)} f(z)\d z.$$
Combining \eqref{MFF}-\eqref{MF}, and noting that
$$r_0:= \min\Big\{-\ff 1 2-\ff d{2k},\ \kk-\ff {dp}{2k'(p-1)}\Big\}>-1,$$
we find a constant $c_2\in (0,\infty)$ such that
\beq\label{NM}    \beg{split} & \d |Y_t^\vv- X_t^\vv|^{2n} \le \d M_t +   c_2 \Big(\big(\E|X_t^\vv-Y_t^\vv|^p\big)^{\ff{2n }p}   + \vv^n  \Big) t^{r_0}\d t\\
 &\quad + c_2  |X_t^\vv-Y_t^\vv|^{2n} \big(t^{r_0} + \scr M\|\si_t\|^2(X_t^\vv)+ \scr M\|\si_t\|^2(Y_t^\vv) \big) \d t,
\ \ \ t\in [0, T].\end{split}\end{equation}
By   the stochastic Gronwall's inequality and Khasminskii's inequality as in the proof of Lemma \ref{L1}, we find a constant $c_3\in (0,\infty)$ such that \eqref{NM} implies
\beq\label{MFN} \beg{split}&\bigg(\E\bigg[\sup_{s\in [0,t]} |Y_s^\vv-X_s^\vv|^n\bigg]\bigg)^2\\
& \le c_3 \int_0^t \Big(\big(\E|X_s^\vv-Y_s^\vv|^p\big)^{\ff{2n }p}   + \vv^n  \Big) s^{r_0}\d s,\ \ t\in [0,T].\end{split}\end{equation}
By \eqref{*0} and $r_0>-1$, the upper bound is finite. So, taking $n=p$ we find a constant $c_4\in (0,\infty)$ such that
$$  \E\bigg[\sup_{t\in [0,T]} |Y_t^\vv-X_t^\vv|^p\bigg] \le c_4\vv^p,\ \ \ \vv\in (0,1].$$
Substituting  into \eqref{MFN} we find a constant $c_5\in (0,\infty)$ such that
\beq\label{MF3}  \E\bigg[\sup_{t\in [0,T]} |Y_t^\vv-X_t^\vv|^n\bigg] \le c_5 \vv^n,\ \ \ \vv\in (0,1].\end{equation}

On the other hand, letting $X_t^{\mu,x}$ solve \eqref{DC} for random  $x$ measurable with respect to $\F_0$, we have
$$Y_t^\vv= X-t^{\mu, X_0^\mu+\vv\eta},\ \ \ X_t= X_t^{\mu, X_0^\mu}.$$
So, by the second inequality in \eqref{DC}, we obtain
$$\E^{\F_0} \bigg[\sup_{t\in [0,T]} |X_t-Y_t^\vv|^n\bigg] = \E^{\F_0} \bigg[\sup_{t\in [0,T]} |X_t^{\mu,x}-X_t^{\mu,y} |^n\bigg]\bigg|_{(x,y)= (X_0^\mu, X_0^\mu+\vv\eta)}
\le C(n) \vv^n |\eta|^n.$$
Noting that $v_t^\vv=\ff{X_t^\vv-X_t}\vv$, this together with \eqref{MF3} implies \eqref{CE} for some $K_n\in (0,\infty).$

(b)  By \eqref{DS} and noting that $X_t=X_t^\mu= X_t^{\mu, X_0^\mu}$,
 we find a constant $C(p,\vv)\in (0,\infty)$ such that
\beq\label{CE'} \E\bigg[\sup_{t\in [0,T]} |v_t^\vv|^p\bigg]= \E\bigg( \E^{\F_0}\Big[\sup_{t\in [0,T]} |v_t^\vv|^p \Big]\bigg)\le C(p,\vv) \E\big[1+|X_0^\mu|^p+|\eta|^p\big]<\infty.\end{equation}
This together with   Lemma \ref{L1} implies  
\beq\label{HT} H_t:= \E \bigg[\sup_{s\in [0,t]} |v_s^\vv-v_s|^p\bigg]<\infty, \ \ \ t\in [0,T].\end{equation}  Let 
\beg{align*}
&\eta_s^\vv:=  \vv\eta + \ff{\si_s(X_s^{\vv})-\si_s(X_s )}\vv- \nn_{v_s^\vv}\si_s(X_s ),\\
&\xi_s^\vv:=  \ff{b_s(X_s^{\vv })-b_s(X_s )}\vv -\nn_{v_s^\vv} b_s(X_s  ) \\
&\qquad + \ff{(h_s*\mu_s)(X_s^{\vv })-(h_s*\mu_s)(X_s )}\vv-
\nn_{v_s^\vv} (h_s*\mu_s) (X_s )\\
& \qquad  + \E\Big[\ff{h_s(z-X_s^{\vv }) - h_s(z-X_s )}\vv- \nn_{v_s^\vv}h_s(z-X_s )\Big]\Big|_{z=X_s^\vv }\\
&\qquad + \E\big[\nn_{v_s^\vv}h_s(z-X_s )- \nn_{v_s^\vv}h_s(z'-X_s )\big]\big|_{z=X_s^\vv, z'=X_s }.
 \end{align*}
Then 
\beg{align*} &v_t^\vv-v_t= \int_0^t \Big(\nn_{v_s^\vv-v_s} b_s(X_s) + \nn_{v_s^\vv-v_s} (h_s*\mu_s) (X_s) +\E\big[\nn_{v_s^\vv-v_s}h_s(z-X_s)\big]\big|_{z=X_s}\Big)\d s\\
&\qquad  + \int_0^t \nn_{v_s^\vv-v_s}\si_s(X_s)\d W_s+ \int_0^t\xi_s^\vv\d s+ \int_0^t \eta_s^\vv\,\d W_s,\ \ \ t\in [0,T].\end{align*}
Let $\theta\in (1, 1+p^{-1})$  such that $\theta p-1<p$.
By $\|\nn b\|_\infty<\infty$ due to $(H)$, \eqref{Y1}, \eqref{X1} for $v_t^\vv-v_t$ in place of $v_s$,
we find a constant $c_1\in (0,\infty)$   and a martingale $M_t$ such that
\beg{align*}&\d |v_t^\vv-v_t|^{\theta p} \le  \d M_t+ c_1\bigg(1+t^{\kk-\ff{dp}{2k'(p-1)}}+\|\nn \si_t(X_t)\|^2  \bigg)|v_t^\vv-v_t|^{\theta p}\d t \\
&+c_1 \big(1+t^{\kk-\ff{dp}{2k'(p-1)}}  \big)\big(\E|v_t^\vv-v_t|^p\big)^\theta \d t +  c_1 \big(|v_t^\vv-v_t|^{\theta p-2}|\eta_t^\vv|^2+|v_t^\vv-v_t|^{\theta p-1} |\xi_t^\vv|\big)\d t,\ \ \ t\in (0,T].\end{align*}
By the stochastic Gronwall inequality and Khasminskii estimate as in the proof of Lemma \ref{L1}, we find  constants $c_2,c_3\in (0,\infty)$ such that
\beg{align*} &\bigg(\E^{\F_0} \bigg[\sup_{s\in [0,t]} |v_s^\vv-v_s|^p\bigg]\bigg)^\theta\\
& \le c_2\E^{\F_0}\int_0^t \Big( \Big(\E\big[|v_s^\vv-v_s|^p\big]\Big)^\theta\big(1+s^{\kk-\ff{dp}{2k'(p-1)}}\big) +
   \big(|v_s^\vv-v_s|^{\theta p-1} |\xi_s^\vv|+|v_s^\vv-v_s|^{\theta p-2}|\eta_s^\vv|^2\big) \Big) \d s \\
&\le \ff 1 2 \bigg(\E^{\F_0} \bigg[\sup_{s\in [0,t]} |v_s^\vv-v_s|^p\bigg] \bigg)^\theta+ c_3 \int_0^t \big(1+s^{\kk-\ff{dp}{2k'(p-1)}}\big)\big(\E\big[|v_s^\vv-v_s|^p\big]\big)^\theta\d s \\
&\qquad + c_3 \E^{\F_0} \bigg[\bigg(\int_0^t |\xi_s^\vv|\d s \bigg)^{\theta p} +  \bigg(  \int_0^t \|\eta_s^\vv\|^2\d s \bigg)^{\ff{\theta p}2}\bigg],\ \ t\in [0,T].\end{align*}
 Combining this with \eqref{HT}, we find a constant $c_4\in (0,\infty)$ such that
\beg{align*} H_t\le  c_4  \bigg(\int_0^t \Big(1+s^{\kk-\ff{dp}{2k'(p-1)}}\Big)  H_s^\theta  \d s \bigg)^{\ff 1 \theta}
  + c_4 (\aa_\vv+\bb_\vv),\ \ \ t\in [0,T],\ \vv\in (0,1],  \end{align*} where
 \beg{align*}  \aa_\vv:=    \E\bigg|  \E^{\F_0}  \bigg(  \int_0^T |\eta_s^\vv|^2\d s \bigg)^{\ff{\theta p}2}\bigg|^{\ff 1 \theta},\ \ \ \ 
 \bb_\vv:= \E    \bigg| \E^{\F_0}\bigg(\int_0^T |\xi_s^\vv|\d s \bigg)^{\theta p} \bigg|^{\ff 1 \theta}.\end{align*}
By \eqref{KL}, this implies $H_T\le c_5(\aa_\vv+\bb_\vv)$ for some constant $c_5\in (0,\infty).$ So, the proof is finished if  we can verify $  \lim_{\vv\downarrow 0}(\aa_\vv+\bb_\vv)=0$.

(c)  Proof of $  \lim_{\vv\downarrow 0}\aa_\vv=0$.  By \eqref{MF} and $(H)$(1), we find a constant $c_1\in (0,\infty)$ such that
\beq\label{ETA} \beg{split} \|\eta_s^\vv\|&\le \Big\|\ff{\si_s(X_s^\vv)- \si_s(X_s)}\vv\Big\|+ |v_s^\vv|\|\nn \si_s(X_s)\| \\
&\le c_1|v_s^\vv| \big(1+ \scr M\|\nn \si_s\| (X_s^\vv) + \scr M\|\nn \si_s\| (X_s) \big).\end{split} \end{equation}
By $(H)$, \eqref{MF}, H\"older's inequality, Krylov's estimate and \eqref{CE},
we find     $c_2 \in (0,\infty)$ such that
\beg{align*}
& \E^{\scr F_0}  \bigg[\bigg(  \int_0^T \|\eta_s^\vv\|^2\d s \bigg)^{\ff{\theta p}2}\bigg]  \\
&\le c_1^{p}  \E^{\F_0} \bigg[\Big(\sup_{s\in [0,T]} |v_s^\vv|^{\theta p}\Big)  \bigg(\int_0^T \big(1+ \scr M\|\nn \si_s\| (X_s^\vv) + \scr M\|\nn \si_s\| (X_s)\big)^2\d s\bigg)^{\ff{\theta p}2}\bigg] \\
& \le c_1^p \bigg(\E^{\F_0}\Big[\sup_{s\in [0,T]} |v_s^\vv|^{2\theta p} \Big]\bigg)^{\ff 1 {2 }} \bigg(\E^{\F_0} \bigg[\bigg(\int_0^T \big(1+ \scr M\|\nn \si_s\| (X_s^\vv) + \scr M\|\nn \si_s\| (X_s)\big)^2\d s\bigg)^{\theta p}\bigg]\bigg)^{\ff 1 {2}} \\
& \le c_2 |\eta|^{\theta p}.\end{align*}
Since $\E[|\eta|^p]<\infty$, by the dominated convergence theorem, $\lim_{\vv\to 0} \aa_\vv =0$ follows if
we can  show that $\P$-a.s.
\beq\label{AF} \lim_{\vv\downarrow 0} \E^{\scr F_0}  \bigg[\bigg(  \int_0^T \|\eta_s^\vv\|^2\d s \bigg)^{\ff{\theta p}2}\bigg]=0.\end{equation}
Since $\nn \si_s(x)$ exists for a.e. $x\in\R^d$, we have
$$\lim_{\vv\downarrow 0} \sup_{|v|\le 1} \Big\|\ff{\si_s(x+\vv v)-\si_s(x)}\vv-\nn_v \si_s(x)\Big\|=0,\ \ \ \text{a.e.}\ x\in\R^d.$$
Since the noise is non-degenerate, for any $s\in (0,T]$, $\L_{X_s}$ is absolutely continuous with respect to the Lebesgue measure, see for instance \cite[Theorem 6.3.1]{BKRS}. So,
$$\lim_{\vv\downarrow 0}|\eta_s^\vv|=0,\ \ \ \P\text{-a.s.},\ s\in   (0,T].$$
By the dominated convergence theorem, \eqref{AF} follows if there exists $\ll>1$ such that $\P$-a.s.
$$  \sup_{\vv\in (0,1]} \E^{\F_0} \bigg[\bigg(\int_0^T |\eta_s^\vv|^{2\ll}\d s\bigg)^{ \theta p }  \bigg]<\infty.$$
By choosing $\ll>1$ such that $(\ll p_i, \ll q_i)\in \scr K$ for $1\le i\le l,$ this follows from $(H)$, \eqref{CE},  \eqref{ETA} and  Krylov's estimate.

(d) Proof of $  \lim_{\vv\downarrow 0}\bb_\vv=0$. By the same reason leading to 
  $\lim_{\vv\to 0} |\eta_s^\vv|=0,$ and noting that $\nn h_t $   is a.e. continuous,  we  have 
\beq\label{XI5} \lim_{\vv\to 0} |\xi_s^\vv|=0,\ \ \ \P\text{-a.s.},\ s\in (0,T]. \end{equation}
By $\|\nn b\|_\infty<\infty$ due to $(H)$, \eqref{X1} and \eqref{CE'}, we find   $c_1,c_2\in (0,\infty)$ such that
\beq\label{CE*} |\xi_s^\vv|\le c_1\Big(|v_s^\vv|+ \big(\E [|v_s^\vv|^p]\big)^{\ff 1 p} \Big)  s^{\kk-\ff{d(p-1)}{2k'(p-1)}}
\le c_2 (1+|v_s^\vv|) s^{\kk-\ff{dp}{2k'(p-1)}}.\end{equation}
By \eqref{KL}, we find $\ll>1$ such that
$$\ll':=\ll \Big(\kk-\ff{dp}{2k'(p-1)}\Big)>-1,$$
so that \eqref{CE} and \eqref{CE*} yield that for some constant $c_3\in (0,\infty)$,
\beg{align*} &\sup_{\vv\in (0,1]} \E^{\F_0} \bigg[\bigg(\int_0^T |\xi_s^\vv|^\ll \d s\bigg)^{\theta p}\bigg]\\
&\le c_2^{\ll\theta p} \bigg(\sup_{\vv\in (0,1]} \E^{\F_0}\Big[\sup_{s\in [0,T]} |v_s^\vv|^{\theta \ll p}\Big]\bigg)\bigg(\int_0^T s^{\ll'}\d s\bigg)^{\theta p} \le c_3 |\eta|^{\ll\theta p}.\end{align*}
Combining this with \eqref{XI5} and the dominated convergence theorem, we drive
$\P$-a.s.
$$\lim_{\vv\to 0}  \bigg| \E^{\F_0}\bigg(\int_0^T |\xi_s^\vv|\d s \bigg)^{\theta p} \bigg|^{\ff 1 \theta}=0$$
and for some constant $c_4\in (0,\infty)$,
$$ \sup_{\vv\in (0,1]}\bigg| \E^{\F_0}\bigg(\int_0^T |\xi_s^\vv|\d s \bigg)^{\theta p} \bigg|^{\ff 1 \theta}\le c_4 |\eta|^p.$$
Since $\E[|\eta|^p]<\infty$, by the dominated convergence theorem again, we obtain
$  \lim_{\vv\downarrow 0}\bb_\vv=0$.
\end{proof}

\section{Proof of Theorem \ref{T1}(2)}

  For any $\eta\in L^p(\OO\to\R^d,\F_0,\P)$, $\mu\in \scr P_k$, consider
$$ \GG_\eta P_t f(\mu):=\lim_{\vv\downarrow 0} \ff{\E[f(X_t^{\vv\eta,\mu})-f(X_t^\mu)]}\vv,\ \ t\in (0,T], f\in \B_b(\R^d).$$
Then 
$$ D_\phi P_t f(\mu)=\GG_{\phi(X_0^\mu)}  P_t f(\mu),\ \ \   \phi\in T_{\mu,p}^*,\ t\in (0,T].$$ 
So, Theorem \ref{T1}(2)  follows from  the following result.

\beg{prp}\label{PP2}  In the situation of Theorem \ref{T1}.  For any $\eta\in L^p(\OO\to\R^d,\F_0,\P)$ and  $\mu\in \scr P_p$,   $ \GG_\eta P_t f(\mu)$ exists and satisfies the following formula for any    $\bb \in C^1([0,t])$ with $\bb_0=0$ and $\bb_t=1:$
\beq\label{BSM2} \beg{split} &  \GG_\eta P_t f(\mu)
= \int_{\R^d\times \R^d} \E\bigg[f(X_t^{\mu,x}) \int_0^t \bb_s'\big\<\zeta_s(X_s^{\mu,x} ) \nn_{v}X_s^{\mu,x},  \d W_s\big\> \bigg]\L_{(X_0^\mu, \eta)}(\d x,\d v)\\
&+  \E\bigg[f(X_t^{\mu}) \int_0^t  \Big\<\zeta_s(X_s^\mu)  \E\big[\<\nn h_s(z-X_s^\mu),\nn_\eta X_s^\mu\>\big]\big|_{z=X_s^\mu},  \d W_s\Big\> \bigg].
\end{split} \end{equation}
Consequently,     there exists a constant $c >0$ such that
\beq\label{GRDI2}\beg{split} & \big| \GG_\eta P_t f(\mu)\big|  \le \ff {c}{\ss t} \big(P_t |f|^{\ff{p}{p-1}}(\mu)\big)^{\ff {p-1}p}(\E[|\eta|^p])^{\ff 1 p},\\
&t\in (0,T], f\in \B_b(\R^d),\mu\in \scr P_p, \eta\in  L^p(\OO\to\R^d,\F_0,\P).\end{split} \end{equation}
\end{prp}

To prove this result, we need the following lemma.

\beg{lem}\label{L21}  Assume $(H)$ and $2\kk-\ff d{k'}>-1.$ Let $p\in [\ff{2k'(\kk+1)}{2k'(\kk+1)-d},\infty)\cap (\ff {k'}{k'-1},\infty)$ and   denote $X_s=X_s^\mu, X_s^\vv= X_s^{\vv\eta,\mu}, s\in [0,T].$ 
\beg{enumerate} \item[$(1)$] For any $\theta\in (1,2\land \ff{k'p}{k'+p})$, there exists  $c(\theta)\in (0,\infty)$ such that
\beq\label{ED1} \sup_{z\in\R^d} \Big(\E\big|h_s(z-X_s)-h_s(z-X_s^{\vv})\big|^\theta \Big)^{\ff 1 \theta}\le \vv c(\theta) s^{\kk-\ff d{2k'}},\ \ \ s\in (0,T],\ \vv\in (0,1].\end{equation}
\item[$(2)$] Let $\mu_t= \L_{X_t}, \mu_t^\vv= \L_{X_t^{\vv}}$ and
 \beq\label{GS}  \Xi_s^\vv:= \zeta_s(X_s^{\vv}) \big\{B_s(X_s^{\vv},\mu_s)- B_s(X_s^{\vv},\mu_s^\vv)\big\},\ \ \ 
  R_t^\vv:= \e^{\int_0^t\<\Xi_s^\vv, \d W_s\>-\ff 1 2\int_0^t |\Xi_s^\vv|^2\d s}.  \end{equation} Then for any $n\in [1,\infty)$,  
\beq\label{LLN} \sup_{t\in [0,T], \vv\in (0,1]} \ff 1 {\vv^n} \E\big[|R_t^\vv-1|^n \big]<\infty. \end{equation}
 \end{enumerate}
\end{lem}

\beg{proof} 

(1) We have
\beg{align*} &\big|h_s(z-X_s)-h_s(z-X_s^{\vv })\big|= \bigg|\int_0^1\ff{\d}{\d r} h_s(z-X_s^{\vv r})\d r\bigg|\\
&=\bigg|\big(\nn_{\nn_\eta X_s^{\vv r}} h\big)(z-X_s^{\vv r})\d r\bigg|.\end{align*}
By H\"older's inequality and \eqref{YB}, we find a constant $c_1\in (0,\infty)$ such that
\beq\label{PP0} \beg{split}&\sup_{z\in\R^d} \Big(\E\big|h_s(z-X_s)-h_s(z-X_s^\vv)\big|^\theta \Big)^{\ff 1 \theta}\\
&\le \bigg(\int_0^1 \E\Big[\big|\nn_\eta X_s^{\vv r}\big|^\theta \big|\nn h_s\big|^\theta(z- X_s^{\vv r})\d r\bigg)^{\ff 1 \theta} \\
&\le \bigg(\int_0^1 \Big(\E\big|\nn_\eta X_s^{\vv r}\big|^{p}\Big)^{\ff 1p}
\Big(\E \big|\nn h_s\big|^{\ff {p\theta}{p-\theta}}(z- X_s^{\vv r})\Big)^{\ff {p-\theta}p} \d r\bigg)^{\ff 1 \theta}\\
&\le c_1 \vv \bigg(\int_0^1 \Big(\E \big|\nn h_s\big|^{\ff {p\theta}{p-\theta}}(z- X_s^{\vv r})\Big)^{\ff {p-\theta}p} \d r\bigg)^{\ff 1 \theta}.\end{split}\end{equation}
By $(H)(3)$ and \eqref{GRD} for $p_1=\ff{k'(p-\theta)}{p\theta}>1$ and $p_2=\infty$, we find a constant $c_2\in (0,\infty)$ such that
\beg{align*} &\Big(\E \big|\nn h_s\big|^{\ff {p\theta}{p-\theta}}(z- X_s^{\vv r}) \Big)^{\ff{p-\theta}p}= \Big(\mu_0^{\vv r}\big(P_s^{\mu^{\vv r}}|\nn h_s|^{\ff {p\theta}{p-\theta}}(z-\cdot)\big)\Big)^{\ff{p-\theta}p}\\
&\le \|\nn h_s\|_{\tt L^{k'}}^{\theta} \|P_s^{\mu^{\vv r}}\|_{\tt L^{\ff{k'(p-\theta)}{p\theta}}\to\tt L^\infty}^{\ff{p-\theta}p}\le c_2 s^{\theta (\kk  -\ff{d }{2k'})},\ \ s\in (0,T],\ r\in (0,1).\end{align*}
Combining this with \eqref{PP0} and  $\theta (\kk  -\ff{d }{2k'})>-1$, we prove  \eqref{ED1} for some     $c(\theta)\in (0,\infty)$. 

  (2) By $(H)$ and \eqref{ED1}, we find a constant $C_1\in (0,\infty)$ such that
\beq\label{XIS'}\beg{split}& |\Xi_s^\vv|\le \|\zeta_s\|_\infty \Big|\E\big[h_s(z-X_s)-h_s(z-X_s^\vv)\big]\big|_{z=X_s}\Big|\\
&\le C_1 \vv s^{\kk-\ff{d}{2k'}},\ \ \ s\in (0,T],\ \vv\in (0,1].\end{split}\end{equation}
Since $2\kk-\ff d{k'}>-1$, for any $n\ge 1,$   there exist $C_2,C_3\in (0,\infty)$ such that
\beg{align*} &\E\bigg|\int_0^t\<\Xi_s^\vv, \d W_s\> \bigg|^{2n} +\bigg(\E \int_0^t |\Xi_s^\vv|^2\d s\bigg)^{2n} \\
&\le C_2\bigg(\E\int_0^t|\Xi_s^\vv|^2\d s\bigg)^n + \bigg(\E \int_0^t |\Xi_s^\vv|^2\d s\bigg)^{2n}\\
&\le C_3 \vv^{2n},\ \ \ t\in [0,T],\ \vv\in (0,1].\end{align*}
Moreover, there exists constants $C_4,C_5\in (0,\infty)$ such that
\beq\label{XIS''}\sup_{t\in [0,T],\vv \in (0,1]}\E \big[|R_t^\vv|^{2n}\big]\le \sup_{t\in [0,T],\vv \in (0,1]} \E\big[\e^{C_4\int_0^t |\Xi_s^\vv|^2\d s}\big]\le C_5.\end{equation}
Combining this with $|\e^r-1|\le (\e^r+1)|r|$ for $r\in R$, we find constants $C_6,C_7\in (0,\infty)$ such that
 \beg{align*} &\E\big[|R_t^\vv-1|^n\big]\le \E\bigg[\big(R_t^\vv+1\big)^n \bigg|\int_0^t\<\Xi_s^\vv,\d W_s\>-\ff 1 2\int_0^t |\Xi_s^\vv|^2\d s\bigg|^n\bigg]\\
 &\le C_6\Big( \E\big[\big(R_t^\vv+1\big)^{2n}\big]\Big)^{\ff 1 2} \bigg(\E  \bigg|\int_0^t\<\Xi_s^\vv,\d W_s\>\bigg|^{2n} + \bigg(\E\int_0^t |\Xi_s^\vv|^2\d s\bigg)^{2n}\bigg)^{\ff 1 2}\\
 &\le C_7 \vv^n,\ \ \ t\in [0,T],\ \vv\in (0,1].\end{align*}
 Therefore, \eqref{LLN} holds.

 \end{proof}

\beg{proof}[Proof of Proposition \ref{PP2}]  By \eqref{DS}, \eqref{Y1}, \eqref{YB} and taking $\bb_s= \ff s t$, we deduce \eqref{GRDI2} from \eqref{BSM2}.  So, it suffices to prove
\eqref{BSM2}.

(a) Let $X_t^{\mu,x}$ solve \eqref{DC}.  
Let
$(P_{s,t}^\mu)_{0\le s\le t\le T}$ be the semigroup associated with    \eqref{DC}, i.e. for  $(X_{s,t}^{\mu,x})_{t\in [s,T]}$ solving \eqref{DC} from time $s$ with $X_{s,s}^{\mu,x}=x$,
\beq\label{TSM} P_{s,t}^\mu f(x):= \E[f(X_{s,t}^{\mu,x})],\ \ t\in [s,T], x\in\R^d.\end{equation}
Then $P_t^\mu= P_{0,t}^\mu$  and
\beq\label{EN0} P_tf(\mu)= \E[f(X_t^\mu)]= \int_{\R^d} P_t^\mu f(x)\mu(\d x),\ \ t\in [0,T], f\in \B_b(\R^d).\end{equation}
 Next, denote $\mu_t=P_t^*\mu=\L_{X_t^\mu}$ and  let ${\bar X}_s^{\vv}$ solve the SDE
 \beq\label{*EN} \d \bar X_s^\vv = \big(b_s(\bar X_s^\vv)+ B_s(\bar X_s^\vv,\mu_s)\big)\d s +\si_s(\bar X_s^\vv)\d W_s,\ \  s\in [0,t], \bar X_0^\vv= X_0^\mu+\vv\eta.\end{equation}
We have
\beg{align*} &\E[ f(\bar X_t^\vv)]= \int_{\R^d} (P_t^\mu)(x) \L_{X_0^\mu+\vv \eta}(\d x)\\
&= \int_{\R^d\times\R^d} P_t^\mu f(x+\vv v)\L_{(X_0^\mu,\eta)}(\d x, \d v),\ \ f\in \B_b(\R^d).\end{align*}
 Combining this with \eqref{EN0} and \eqref{B0}, and applying the dominated convergence theorem, we obtain
 \beq\label{EN2}\beg{split} & \lim_{\vv\to 0}  \ff{\E[f(\bar X_t^\vv)] - P_t f(\mu)}\vv =  \int_{\R^d\times\R^d} \nn_{v} P_t^\mu f(x)\L_{(X_0^\mu,\eta)}(\d x, \d v)\\
 &= \int_{\R^d\times\R^d} \E\bigg[f(X_t^{\mu,x}) \int_0^t \bb_s' \<\zeta_s(X_s^{\mu,x}) \nn_{v} X_s^{\mu,x}, \d W_s\>\bigg]\L_{(X_0^\mu,\eta)}(\d x, \d v).\end{split}\end{equation}
By \eqref{XIS'} and $2\kk-\ff d{k'}>-1$,  we find a constant $c_2\in (0,\infty)$ such that 
 $$\int_0^T |\Xi_s^\vv|^2\d s\le c_2.$$ So,
 by  Girsanov's theorem, $\Q_t^\vv:= R_t^\vv\P$ is a probability measure under which
$$\tt W_r^\vv:=W_r- \int_0^r\zeta_s(X_s^\vv) \big(B_s(X_s^\vv,\mu_s)- B_s(X_s^\vv,\mu_s^\vv)\big)\d s,\ \ r\in [0,t]$$ is a Brownian motion.
 Reformulate the SDE  for $X_s^\vv$ as
$$ \d X_s^\vv = \big(b_s(X_s^\vv)+B_s(X_s^\vv, \mu_s) \big)\d s+\si_s(X_s^\vv)\d \tt W_s^\vv,\ \ X_0^\vv= \bar X_0^\vv.$$
By the weak well-posedness we obtain
$$ \E[f(\bar X_t^\vv)]= \E[R_t^\vv f(X_t^\vv)],\ \ f\in \B_b(\R^d).$$
Thus,
 \beg{align*} & \ff{\E[ f(X_t^\vv)]- \E[f(\bar X_t^\vv)]}\vv =  \ff{\E[f(X_t^\vv)(1-R_t^\vv)]}\vv= I_1(\vv)+ I_2(\vv),\\
& I_1(\vv):= \E\bigg[f(X_t^\mu) \ff{1-R_t^\vv}\vv\bigg],\ \ I_2(\vv):= \E\bigg[\big(f(X_t^\vv)-f(X_t^\mu)\big) \ff{1-R_t^\vv}\vv\bigg].\end{align*}
So, it remains to verify
\beq\label{I1} \lim_{\vv\to 0} I_1(\vv)
 = \E\bigg[ f(X_t^\mu) \int_0^t \big\<\zeta_s(X_s)  \E[\<\nn h_s(z-X_s),\nn_\eta X_s^\mu\>]|_{z=X_s},  \d W_s\big\>\bigg],\end{equation}
\beq\label{I2} \lim_{\vv\to 0} I_2(\vv)=0,\end{equation}
By the definition of $R_t^\vv$ in \eqref{GS} and the estimate \eqref{LLN},  we may apply the dominated convergence theorem to derive
$$\lim_{\vv\downarrow 0} \ff{R_t^\vv-1}\vv=  \int_0^t \big\<\zeta_s(X_s)  \E[\<\nn h_s(z-X_s),\nn_\eta X_s\>]|_{z=X_s},  \d W_s\big\>$$
in $L^1(\P)$ (indeed in $L^n(\P)$ for any $n\in [1,\infty)$).  So, \eqref{I1} holds.

(b) Proof of \eqref{I2}. When $f\in C_b(\R^d)$, it follows from \eqref{LLN}  and the dominated convergence theorem. Below we prove
it for $f\in \B_b(\R^d)$ by using a time-shift argument as in the proof of \cite[Proposition 4.3]{W23}.

By $2\kk-\ff d{k'}>-1$,  \eqref{XIS'} and \eqref{XIS''},   we find constants $K_1,K_2\in (0,\infty)$ such that for any
$T\ge t\ge r\ge 0$,
\beg{align*} &\ff 1 \vv \E\big[|R_t^\vv-R_r^\vv||\big]= \ff 1 \vv \E\big[|R_r^\vv(\e^{\int_r^t\<\Xi_s^\vv,\d W_s\>-\ff 1 2\int_r^t |\Xi_s^\vv|^2\d s}-1)|\big]\\
&\le \ff 1\vv \Big(\E\big[|R_r^\vv|^2\big]\Big)^{\ff 1 2}\Big(\E\big[\e^{2\int_r^t\<\Xi_s^\vv,\d W_s\>- \int_r^t |\Xi_s^\vv|^2\d s}-1\big]\Big)^{\ff 1 2}\\
&\le \ff{K_1}\vv \Big(\e^{C_1^2\vv^2\int_r^t s^{2\kk-\ff d{k'}}\d s}-1\Big)^{\ff 1 2}\le K_2 \ss{t^{\kk+1-\ff d{k'}}- r^{\kk+1-\ff d{k'}}},\ \ \vv\in [0,1].\end{align*}
Therefore,
\beq\label{LN11} \lim_{r\uparrow t} \sup_{\vv\in (0,1]} \E\bigg[\ff{|R_t^\vv-R_r^\vv|}\vv\bigg]=0,\ \ t\in (0,T].\end{equation}
By the same reason leading to \eqref{GRD}, we find a constant $c\in (0,\infty)$ such that
\beq\label{GRD'} \|\nn P_{r,t}^\mu f\|_\infty\le c(t-r)^{-\ff 1 2}\|f\|_\infty,\ \ \ T\ge t>r\ge 0,\ f\in \B_b(\R^d),\ \mu\in \scr P.\end{equation}
Then by the Markov property,
\beq\label{NL0} \beg{split}& |\E[f(X_t^\vv)-f(X_t^\mu)|\F_{r}]|= |(P_{r,t}^{\mu^\vv} f)(X_{r}^\vv)- (P_{r, t}^\mu f)(X_{r}^\mu)| \\
&\le  |(P_{r,t}^{\mu^\vv} f)(X_{r}^\vv)- (P_{r, t}^{\mu^\vv} f)(X_{r}^\mu)|+  |(P_{r,t}^{\mu^\vv} f)(X_{r}^\mu)- (P_{r, t}^\mu f)(X_{r}^\mu)|\\
&\le  c \|f\|_\infty\E\bigg[\ff{|X_{r}^\vv- X_{r}^\mu|}{\ss{t-s}} \bigg] +  |(P_{r,t}^{\mu^\vv} f)(X_{r}^\mu)- (P_{r, t}^\mu f)(X_{r}^\mu)|.\end{split}\end{equation}
On the other hand, let $(\tt X_{r,s}^\vv)_{s\in [r,t]}$ solve the SDE
$$\d \tt X_{r,s}^\vv = b_s(\tt X_{r,s}^\vv, \mu_s^\vv)\d s +\si_s(\tt X_{r,s}^\vv)\d W_s,\ \ \tt X_{r,r}^\vv= X_{r}^\mu, s\in [r,t].$$
We have
$$P_{r,t}^{\mu^\vv} f(X_{r}^\mu)= \E\big[f(\tt X_{r,t}^\vv)\big|\F_{r}\big],\ \ P_{r,t}^\mu f(X_{r}^\mu)= \E\big[f(X_t^\mu)|\F_{r}\big].$$
By \eqref{XIS'} and  Girsanov's theorem,
$$R_{r,t}^\vv:=\e^{\int_{r}^t\<\Xi_s^\vv,\d W_s\>-\ff 1 2\int_{r}^t |\Xi_s^\vv|^2\d s},\ \ \ t\in [r,T]$$
is a martingale, and
$$\tt W_s:= W_s- \int_{r}^s \Xi_\theta^\vv \d \theta,\ \ s\in [r,t]$$
is a Brownian motion under $\Q_{r,t}:= R_{r,t}^\vv \P$.
Reformulating the SDE for $(X_{s}^\mu)_{s\in [r,t]}$ as
$$\d X_s^\mu= \big\{b_s(X_s^\mu)+B_s(X_s^\mu, \mu^\vv)\d s +\si_s(X_s^\mu) \d \tt W_s,\ \ X_{r}^\mu= \tt X_{r,r}^\vv,\ \ s\in [r,t], $$
by the weak uniqueness we obtain
$$P_{r,t}^{\mu^\vv} f(X_{r}^\mu)= \E\big[R_{r,t}^\vv f(X_t^\mu)\big|\F_{r}\big],$$
so that by   Pinsker's inequality, $2\kk-\ff d{k'}>-1$  and \eqref{XIS'}, we find a  constant  $c_2>0$ such that
\beq\label{NL*} \beg{split} &|(P_{r,t}^{\mu^\vv} f)(X_{r}^\mu)- (P_{r, t}^\mu f)(X_{r}^\mu)|^2\le \|f\|_\infty \big|\E[|1-R_{r,t}^\vv|\big|\F_{r}]\big|^2\\
&\le 2\|f\|_\infty\E_{\Q_{r,t}}\big[ \log R_{r,t}^\vv\big|\F_0\big]
= \|f\|_\infty\int_{r}^t \E_{\Q_{r,t}}\big[|\Xi_s^\vv|^2\big| \F_{r}\big] \d s\\
&\le   c_2 \|f\|_\infty\vv^2,\ \ \vv\in (0,1],\ 0<r<t\le T.\end{split}\end{equation}
Combining this with Lemma \ref{L2},   \eqref{LLN} and \eqref{NL0},   we find constants $c_3,c_4>0$ such that
\beg{align*}&  \bigg|\E\Big[\big(f(X_t^\vv)-f(X_t)\big)\ff{1-R_r^\vv}\vv\Big] \bigg|\\
&\le \bigg(\E\Big|\E\big[f(X_t^\vv)-f(X_t)\big|\F_r\big]\Big|^2\bigg)^{\ff 1 2} \bigg(\E\Big[\ff{|1-R_r^\vv|^2}{\vv^2} \Big]\bigg)^{\ff 1 2}\\
&\le  c_4\|f\|_\infty\bigg(  \ff{\E[|X_{r}^\vv-X_{r}^\mu|]}{\ss{t-r}}  \bigg)^{\ff 1 2}   + c_4\|f\|_\infty\vv\\
 & \le c_5\|f\|_\infty\Big(\ff{\vv}{t-r}\Big)^{\ff 1 2},\ \ \vv\in (0,1], t\in [0,T].\end{align*}
This together  with \eqref{LN11} yields
$$\lim_{\vv\downarrow  0} I_2(\vv)\le \lim_{r\uparrow  t} \lim_{\vv\downarrow 0} \Bigg(\bigg|\E\Big[\big(f(X_t^\vv)-f(X_t)\big)\ff{1-R_r^\vv}\vv\Big]\bigg|+2 \|f\|_\infty \E\Big[\ff{|R_t^\vv-R_r^\vv|}\vv\Big]\Bigg) =0. $$
  \end{proof}

\end{document}